\documentclass[12pt]{amsart}

\usepackage{times}

\author{Piotr \'Sniady}
\address{Institute of Mathematics,
University of Wroclaw, pl.\ Grunwaldzki 2/4, 50-384 Wroclaw, Poland}
\email{Piotr.Sniady@math.uni.wroc.pl}

\title[Inequality for Voiculescu's free entropy]{Inequality for Voiculescu's
free entropy in terms of Brown measure}

\theoremstyle{plain}
\newtheorem{lemma}{Lemma}
\newtheorem{theorem}[lemma]{Theorem}
\newtheorem{proposition}[lemma]{Proposition}

\theoremstyle{definition}

\theoremstyle{remark}

\newcommand{\A}{{\mathcal{A}}}
\newcommand{\E}{{\mathbb{E}}}
\newcommand{\El}{{\mathcal{L}}}

\newcommand{\C}{{\mathbb{C}}}
\newcommand{\R}{{\mathbb{R}}}
\newcommand{\M}{{\mathcal{M}}}

\newcommand{\N}{{\mathbb{N}}}
\newcommand{\gwia}{^{\star}}
\newcommand{\random}{\M_N\big(\El^{\infty-}(\Omega) \big)}
\newcommand{\dd}{\ {\mathrm{d}}}

\DeclareMathOperator{\tr}{tr}
\DeclareMathOperator{\Tr}{Tr}
\DeclareMathOperator{\od}{od}

\DeclareMathOperator{\Var}{Var}
\DeclareMathOperator{\vol}{vol}

\newcommand{\sa}{{\operatorname{sa}}}

\newcommand{\dia}{{\operatorname{d}}}
\newcommand{\su}{{\operatorname{sut}}}
\newcommand{\sut}{{\operatorname{sut}}}
\newcommand{\uu}{{\operatorname{U}}}

\begin{document}

\bibliographystyle{alpha}

\begin{abstract}
We study Voiculescu's microstate free entropy for a single non--selfadjoint
random variable. The main result is that certain additional constraints on
eigenvalues of microstates do not change the free entropy. Our tool is
the method of random regularization of Brown measure which was studied recently
by Haagerup and the author. As a simple application we present an upper
bound for the free entropy of a single non--selfadjoint operator in terms of
its Brown measure and the second moment. We furthermore show that this
inequality becomes an equality for a class of $DT$--operators which was
introduced recently by Dykema and Haagerup.
\end{abstract}

\maketitle

\section{Introduction}
The microstate free entropy $\chi$ was introduced by Voiculescu
\cite{VoiculescuPart2} as a tool for the study of some non--commutative
systems. Roughly speaking, it answers the question how many finite matrices
have nearly the same moments as a given non--commutative random variable. It
has turned out to have very powerful applications (cf.\
\cite{Ge1997,VoiculescuPart3}), however it is not an easy object to deal with.

One of the reasons of these difficulties is that currently there are no
general methods for the computation of the free entropy in concrete cases.
Exact formulas were found for the free entropy of a single selfadjoint random
variable, for tuples of free random variables \cite{VoiculescuPart2} and for
$R$--diagonal elements \cite{NicaShlyakhtenkoSpeicher1999}.

In this article we present a method which hopefully will be
useful for calculating and estimating free entropy in many concrete cases.
The main idea is to change the definition of microstates $\Gamma$ which
approximate a single non--selfadjoint random variable $x$ in such a way that
it does not change the value of the free entropy $\chi(x)$.

The original sets $\Gamma$ consisted of all matrices which---informally
speaking---had almost the same moments as a given random variable $x$, while
our new sets $\tilde{\Gamma}$ will consist of these matrices in $\Gamma$ which
additionally have similar eigenvalues to the Brown measure of $x$. In order to
show that $\Gamma$ and $\tilde{\Gamma}$ give rise to the same free entropy we
use the method of random regularization of Brown spectral measure which was
introduced by Haagerup \cite{Haagerup2001} and further developed by the author
\cite{Sniady2001}.

As an application we present a new upper bound for the free entropy of a single
random variable $x$ in terms of its Brown measure and second moment. We
show also that for a class of $DT$--operators which was introduced recently by
Dykema and Haagerup \cite{DykemaHaagerup2001} this inequality becomes equality.

\section{Preliminaries}
\label{sec:preliminaries}
\subsection{Non--commutative probability spaces}
A non--commutative probability space is a pair $(\A,\phi)$, where $\A$ is a
$\star$--algebra and $\phi$ is a normal, faithful, tracial state on $\A$.
Elements of $\A$ will be referred to as non--commutative random variables and
state $\phi$ as expectation value.

One of the simplest examples is the set $\M_N$ of all complex--valued $N\times
N$ matrices equipped with a normalized trace $\tr$ given by
$\tr m=\frac{1}{N} \Tr m$, where $m\in\M_N$ and $\Tr$ is the usual trace.

\subsection{Microstate free entropy.}
The original definition of Voiculescu's free entropy $\chi^\sa(x_1,\dots,x_n)$
allowed to compute the free entropy of a tuple of non--commutative
self--adjoint random variables. The considered in this article free entropy
$\chi(x)$ of a non--selfadjoint random variable is connected with the original
definition by $$\chi(x)=\chi^\sa(\Re x,\Im x).$$

Let $x$ be a non--commutative random variable,
 $\epsilon>0$, $R>0$ be real numbers and $k>0$ be integer.
We define the sets \cite{VoiculescuPart2,NicaShlyakhtenkoSpeicher1999}
\begin{multline}
  \Gamma_R(x;k,N,\epsilon)=\Big\{m\in \M_N: \|m\|\leq R \text{ and} \\
   |\tr(m^{s_1}\cdots m^{s_p})-
  \tau(x^{s_1} \cdots x^{s_p})|<\epsilon\\  \text{for all } p\leq k \text{ and
} s_1,s_2,\dots, s_p\in \{1,\star\} \Big\}. \label{eq:defgamma}
\end{multline}
Define next
\begin{equation}
\chi_R(x;k,\epsilon)=\limsup_{N\to\infty}
\left[\frac{1}{N^2}\log\vol \Gamma_R(x;k,N,\epsilon)
+ \log N\right],
\label{eq:defchir}
\end{equation}
where $\vol$ is a Lebesgue measure on $\M_N$ as described in
(\ref{eq:defvol}).
Lastly, the free entropy is defined by
\begin{equation}
\chi(x)=\sup_R \inf_{k,\epsilon} \chi_R(x;k,\epsilon).
\label{eq:defchi}
\end{equation}

Since $\chi_R(x;k,\epsilon)$ is a decreasing function of $k$ and an increasing
function of $\epsilon$, hence we have the following simple lemma.
\begin{lemma}
\label{lem:pierwszy}
Let a non--commutative random variable $x$ and a number $R>0$ be given. Then
there exists a sequence $(\epsilon_N)$ of non--negative numbers and a sequence
$(k_N)$ of natural numbers such that
$\lim_{N\rightarrow\infty} \epsilon_N =0$,
$\lim_{N\rightarrow\infty} k_N=\infty$ and
\begin{equation}
\chi_R(x)
\leq\limsup_{N\rightarrow\infty} \chi_R(x;k_N,N,\epsilon_N).
\label{eq:nawygnaniu}
\end{equation}
\end{lemma}

\subsection{Fuglede--Kadison determinant and Brown measure}
Let a non--commutative probability space $(\A,\phi)$ be given.
For $x\in\A$ we define its Fuglede--Kadison determinant $\Delta(x)$ by
\cite{FugledeKadison}
$$\Delta(x)=\exp\left[ \phi( \ln |x| )\right]$$
and its Brown measure \cite{Brown} to be the
Schwartz distribution on $\C$ given by
$$\mu_x= \frac{1}{2\pi} \left( \frac{\partial^2}{\partial (\Re\lambda)^2}
+\frac{\partial^2}{\partial (\Im\lambda)^2} \right) \ln \Delta (x-\lambda).$$
One can show that in fact $\mu_x$ is a positive probability measure on $\C$.

\begin{lemma}
The Brown measure of a matrix $m\in\M_N$ with respect to the state $\tr$
is a probability counting measure on the set of eigenvalues of $m$:
$$\mu_m=\frac{1}{N} \sum_{i=1}^N \delta_{\lambda_i}, $$
where $\lambda_1,\dots,\lambda_N$ are the eigenvalues of $m$ counted with
multiples.
\end{lemma}

In the following we will be interested in studying the random measure
$\omega\mapsto \mu_{A(\omega)}$ for a random matrix $A\in\random$.
This random measure is called the empirical eigenvalues distribution.

\subsection{Convergence of $\star$--moments}
Let a sequence $(A_N)$ of random matrices (where $A_N\in\random$), a
non--commutative probability space $(\A,\phi)$ and $x\in\A$ be given.
We say that the sequence $A_N$ converges to $x$ in $\star$--moments almost
surely if for every $n\in\N$ and $s_1,\dots, s_n\in\{1,\star\}$ we have that
$$\lim_{N\rightarrow\infty} \tr_N[ A_N^{s_1} \cdots
A_N^{s_n} ] = \phi( x^{s_1} \cdots x^{s_n} )$$
holds almost surely.

\subsection{Random regularization of Brown measure}
We say that a random matrix
$$G_N=(G_{N,ij})_{1\leq i,j\leq N}\in \random$$
is a standard Gaussian random matrix if
$$\big(\Re G_{N,ij}\big)_{1\leq i,j\leq N},
\big(\Im G_{N,ij}\big)_{1\leq i,j \leq N}$$
are independent Gaussian variables with mean zero and variance $\frac{1}{2 N}$.

\begin{theorem}
\label{theo:regularyzacja}
Let $(A_N)$ be a sequence of random matrices, $A_N\in\random$, which
converges in $\star$--moments to $x$ almost surely.

Let furthermore $(G_N)$ be a
sequence of independent standard Gaussian matrices which is independent of
$(A_N)$.
There exists a sequence $(t_N)$ of real numbers such that
$\lim_{N\rightarrow\infty} t_N=0$ and such that the sequence of empirical
eigenvalues distributions  $\mu_{A_N+t_N G_N}$ converges weakly to $\mu_x$
almost surely.

There also exists a sequence $(B_N)$ of non--random matrices $B_N\in\M_N$ such
that $\lim_{N\rightarrow\infty} \|B_N\|=0$ and such that the sequence of
empirical eigenvalues distributions $\mu_{A_N+B_N}$ converges weakly to $\mu_x$
almost surely.
\end{theorem}
\begin{proof}
The first part was was proved in \cite{Sniady2001}.

For the second part of the theorem let us define $B_N=t_N G_N(\omega)$. Since
$\limsup_{N\rightarrow\infty} \|G_N\|<\infty$ holds almost surely \cite{Geman},
hence so defined sequence $B_N$ fulfills the hypothesis of the
theorem almost surely. \end{proof}

\section{The main result: improved microstates $\tilde{\Gamma}$}
\label{sec:themain}
Let $x$ be a non--commutative random variable,
 $\epsilon>0$, $R>0$ be real numbers and $k>0$, $m>0$ be integers.
In the full analogy with (\ref{eq:defgamma})---(\ref{eq:defchi}) we define
improved microstates $\tilde{\Gamma}$ and improved free entropy~$\tilde{\chi}$:
\begin{multline}
  \tilde{\Gamma}_R(x;k,N,\epsilon,l,\theta)=\bigg\{m\in
\Gamma_R(x;k,N,\epsilon):\\
\left| \int_\C z^i \bar{z}^j d\mu_m - \int_\C z^i \bar{z}^j \dd\mu_x\right|
<\theta \quad\text{for } i,j\leq l
\bigg\}, \label{eq:defgammatilde}
\end{multline}
\begin{equation}
\tilde{\chi}_R(x;k,\epsilon,l,\theta)=\limsup_{N\to\infty}
\left[\frac{1}{N^2}\log\vol \tilde{\Gamma}_R(x;k,N,\epsilon,l,\theta)
+ \log N\right],
\label{eq:defchirtilde}
\end{equation}
\begin{equation}
\tilde{\chi}(x)=\sup_R \inf_{k,\epsilon,l,\theta}
\tilde{\chi}_R(x;k,\epsilon,l,\theta). \label{eq:defchitilde}
\end{equation}

\begin{theorem}
\label{theo:glowne}
For every non--commutative random variable $x$ we have
$$\chi(x)=\tilde{\chi}(x).$$
\end{theorem}
\begin{proof}
Since $\tilde{\Gamma}_R(x; k,N,\epsilon,l,\theta)\subseteq \Gamma_R(x;
k,N,\epsilon)$ hence $\tilde{\chi}(x)\leq \chi(x)$ follows easily.

Let $(\epsilon_N)$ and $(k_N)$ be the sequences given by Lemma
\ref{lem:pierwszy}. Let $(A_N)$ be a sequence of independent random matrices
such that the distribution of $A_N$ is the uniform distribution on the set
$\Gamma_R(x; k_N,N,\epsilon_N)$ and let $(B_N)$ be the sequence given by
Theorem \ref{theo:regularyzacja}.

Since $\|B_N\|$ converges to zero, hence there exists $R'>0$, a sequence of
positive numbers $(\epsilon'_N)$ which converges to zero and a sequence of
natural numbers $(k'_N)$ which diverges to infinity such that
$$\Gamma_R(x;k_N,N,\epsilon_N)+B_N\subseteq\Gamma_{R'}(x;k'_N,N,\epsilon_N')$$
holds for every $N\in\N$,
where $\Gamma_R(x;k_N,N,\epsilon_N)+B_N$ denotes a translation of the set
$\Gamma_R(x;k_N,N,\epsilon_N)$ by the vector $B_N$.

Since random measures $\omega\mapsto\mu_{A_N(\omega)+B_N}$ converge weakly to
$\mu_x$ in probability, hence for any $\theta>0$ and integer $l>0$ we have that
$$\lim_{N\rightarrow\infty}
P\big(\omega:A_N(\omega)+B_N\not\in\tilde{\Gamma}_{R'}(x;
k'_N,N,\epsilon_N',l,\theta) \big) =0.$$
Since the Lebesgue measure is translation--invariant, hence for any $\theta>0$
and integer $l>0$ we have
$$\limsup_{N\rightarrow\infty} \frac{\vol \Gamma_R(x;
k_N,N,\epsilon_N)}{\vol \tilde{\Gamma}_{R'}(x; k'_N,N,\epsilon'_N,l,\theta)}
\leq 1, $$
or equivalently
\begin{equation} \label{eq:christmasa} \limsup_{N\rightarrow\infty}  \big(
\chi_R(x; k_N,N,\epsilon_N) - \tilde{\chi}_{R'}(x; k'_N,N,\epsilon'_N,l,\theta)
\big)\leq 0.
\end{equation}

For any $\epsilon>0$ and integer $k>0$ there exists $N_0$ such that for any
$N>N_0$ we have $\epsilon_N<\epsilon$ and $k_N>k$, hence for $N>N_0$
\begin{equation}
\label{eq:christmasb}
\tilde{\chi}_{R'}(x; k'_N,N,\epsilon'_N,l,\theta) \leq
\tilde{\chi}_{R'}(x; k,N,\epsilon,l,\theta).
\end{equation}

Inequalities \eqref{eq:nawygnaniu}, \eqref{eq:christmasa} and
\eqref{eq:christmasb} combine to give
$$\chi(x)\leq \tilde{\chi}(x).$$
\end{proof}

\section{Application: upper bound for free entropy}
In this section we present a new inequality for the free entropy of a single
non--selfadjoint random variable. The main idea is to write matrices from
microstates $\tilde{\Gamma}$ in the upper triangular form and then to find
constraints on diagonal and offdiagonal entries.

\subsection{Pull--back of the Lebesgue measure on $\M_N$}
We denote by
$\M_N^\dia=\{m\in\M_N: m_{ij}=0 \mbox{ if } i\neq j\}$ the set of diagonal
matrices and by $\M_N^\su=\{m\in\M_N: m_{ij}=0 \mbox{ if } i\geq j\}$ the set
of all strictly upper triangular matrices.

We can regard $\M_N$ and $\M_N^\su$ as real Euclidean spaces with a scalar
product $\langle x,y\rangle=\Re \Tr x y\gwia$ and thus equip them
with Lebesgue measures
\begin{equation}
\vol=\prod_{1\leq i,j\leq N} \dd \Re m_{ij} \dd \Im m_{ij}
\label{eq:defvol}
\end{equation}
and
$$\vol^\su= \prod_{1\leq i<j\leq N} \dd \Re m_{ij} \dd \Im m_{ij}$$
respectively.

We have a clear isomorphism $\M_N^\dia=\{ (\lambda_1,\dots,\lambda_N)\in \C^N \}$
and we equip it with a measure
$$\vol^\dia=\frac{ \pi^{\frac{N^2-N}{2}} }{\prod_{1\leq i\leq N} i!}
 \prod_{1\leq i<j\leq N} |\lambda_i-\lambda_j|^2 \prod_{1\leq i\leq N} \dd\Re
\lambda_i \dd\Im \lambda_i. $$

We also denote by $U_N$ the set of unitary $N\times N$ matrices equipped with
the Haar measure $\vol^\uu$ normalized in such a way that it is a probability
measure.

\begin{proposition}
\label{prop:pullback}
For every $N$ the measure $\vol^\dia\times  \vol^\su \times \vol^\uu$ is a
pull--back of the measure $\vol$ with respect to the map
$$\M_N^\dia \times \M_N^\su \times U_N
\ni (d, m, u) \mapsto u(d+m)u^{-1} \in \M_N. $$
\end{proposition}
\begin{proof}
This result is due to Dyson and can be extracted from Appendix A.35 of
\cite{Mehta}. \end{proof}

\subsection{Diagonal entropy $\hat{\chi}^\dia$}
Let $x$ be a non--commutative random variable. In the
following we define an auxiliary quantity $\hat{\chi}^\dia(\nu)$ which would
answer the question how many diagonal matrices (with respect to the measure
$\vol^\dia$) have almost the same Brown measure as $x$.
\begin{multline*}
\hat{\Gamma}^\dia_R (x;N,l,\theta)=
\bigg\{m \in \M_N^\dia: \|m\|\leq R \text{ and} \\ \bigg| \int_\C z^i
\bar{z}^j\dd\mu_m - \int_\C z^i \bar{z}^j\dd\mu_x \bigg|< \theta \text{ for
all } i,j\leq l \bigg\}; \end{multline*}
$$\hat{\chi}_R^\dia(x;l,\theta)=\lim_{N\to\infty}
\left[\frac{1}{N^2}\log \vol^\dia \hat{\Gamma}_R^\dia(x;l,\theta) +
\frac{\log N}{2} \right],$$
$$\hat{\chi}^\dia(x)=\sup_R \inf_{l,\theta} \hat{\chi}^\dia(x;
l,\theta).$$

\begin{theorem}
\label{theo:diagonalne}
For any non--commutative random variable $x$ we have
$$\hat{\chi}^\dia(x)=\int_{\C} \int_{\C} \log |z_1-z_2| \dd\mu_x(z_1)
\dd\mu_x(z_2)+\frac{3}{4}+\frac{\ln \pi}{2}. $$ \end{theorem}
\begin{proof}
Proof follows exactly the proof of Proposition 4.5 of \cite{VoiculescuPart2},
but since we are dealing with measures on $\C$ and the original proof concerns
measures on $\R$, we have to replace Lemma 4.3 in \cite{VoiculescuPart2} by
Theorem 2.1 of \cite{Hadwin}.
\end{proof}

\subsection{Offdiagonality}
If $x$ is a non--commutative random variable we define its offdiagonality
$\od_x$ by $$\od_x=\tau(xx\gwia)-\int_\C |z|^2 d\mu_x(z).$$
This quantity can be regarded as a kind of a non--commutative variance. Since
$\od_x=0$ if and only if $x$ is normal, hence offdiagonality of $x$ can be also
regarded as a kind of a distance of $x$ to normal operators.

For an upper--triangular matrix $m\in\M_N(\C)$ ($m_{ij}=0$ if $i>j$) its
offdiagonality is equal to the (normalized) sum of squares of the offdiagonal
entries:
$$\od_m=\frac{1}{N} \sum_{1\leq i<j\leq N} |m_{ij}|^2.$$

\begin{proposition}
\label{prop:pozadiagonalne}
For any $o>0$ we have that
$$\lim_{N\rightarrow\infty} \left[ \frac{1}{N^2} \log \vol^\sut
\{m\in\M_N^\sut: \od_m\leq o\}+\frac{\log N}{2} \right]=\frac{1}{2}+\frac{\log
2\pi o}{2}. $$ \end{proposition}
\begin{proof}
It is enough to notice that
$\{m\in\M_N^\sut: \tr mm\gwia\leq o \}$ is a
$N(N-1)$--dimensional ball with radius $\sqrt{oN}$, hence its volume is equal
to $$\pi^{\frac{N(N-1)}{2}} \left[ \Gamma\left( \frac{N (N-1)}{2}+1
\right)\right]^{-1} (oN)^{\frac{N (N-1)}{2}}.$$
\end{proof}

\subsection{The main inequality}
\begin{theorem}
Let $x$ be a non--commutative random variable. Then
\begin{equation}
\label{eq:glowne}
\chi(x)\leq \int_{\C} \int_{\C} \log |z_1-z_2| \dd\mu_x(z_1)
\dd\mu_x(z_2)+\frac{5}{4}+\ln \pi\sqrt{2 \od_x}.
\end{equation}
\end{theorem}
\begin{proof}
Proof is a direct consequence of Theorem \ref{theo:glowne}, Proposition
\ref{prop:pullback}, Theorem \ref{theo:diagonalne} and Proposition
\ref{prop:pozadiagonalne}.
\end{proof}

\subsection{Free entropy of $DT$--operators}
For any compactly supported probability measure $\nu$ on $\C$ and $o\geq 0$
Dykema and Haagerup consider an operator $x$ which is said to be $DT(\nu,o)$
\cite{DykemaHaagerup2001}. This operator is implicitly defined to be the
expected $\star$--moment limit of random matrices
\begin{equation}
A_N=D_N+\sqrt{o} T_N,
\label{eq:definicjaan}
\end{equation}
where $D_N$ is a diagonal random matrix with eigenvalues $\lambda_1,\dots,
\lambda_n$ which are i.i.d.\ random variables with distribution given
by $\nu$ and
\begin{equation}
T_N=\left[ \begin{array}{ccccc}
0 & g_{1,2} &   \cdots & g_{1,n-1} & g_{1,n} \\
0 & 0 & \cdots & g_{2,n-1} & g_{2,n} \\
\vdots&                 &       \ddots & \vdots   & \vdots \\
& & & 0 & g_{n-1,n} \\
0&       & \cdots &             0               & 0
\end{array} \right],
\label{eq:utm}
\end{equation}
is an upper--triangular random matrix where $(\Re g_{i,j},\Im
g_{i,j})_{1\leq i<j \leq N}$ are i.i.d.\ $N\left(0,\frac{1}{N}\right)$ random
variables. We recall that the Brown measure and offdiagonality of $x$ are given
by $\mu_x=\nu$ and $\od_x=o$.

\begin{theorem}
For any compactly supporded probability measure $\nu$ on $\C$ and any number
$o>0$ if $x$ is a $DT(\nu,o)$ then the inequality (\ref{eq:glowne}) becomes
equality.
\end{theorem}
\begin{proof}
Let us fix $R>0$. Similarly as in Lemma \ref{lem:pierwszy} let $(\theta_N)$ be
a sequence of non--negative numbers and $(l_N)$ be a sequence of natural
numbers such that $$\lim_{N\rightarrow\infty} \theta_N =0,
\qquad \lim_{N\rightarrow\infty} l_N=\infty,$$
$$\hat{\chi}^\dia_R(x) \leq\limsup_{N\rightarrow\infty}
\hat{\chi}^\dia_R(x;N,l_N,\theta_N).$$

In definition (\ref{eq:definicjaan}) of $A_N$ let us change $D_N$ to be
any (non--random) element of the set $\hat{\Gamma}^\dia_R(x; N, l_N,\theta_N)$.
From results of Dykema and Haagerup \cite{DykemaHaagerup2000} it follows that
despite this change the sequence $A_N$ still converges in expected
$\star$--moments to $x$:
\begin{equation} \lim_{N\rightarrow\infty} \E \tr (A_N^{s_1} \cdots
A_N^{s_k})= \phi(x^{s_1} \cdots x^{s_k})
\label{eq:zbieznosca}
\end{equation}
for any $k\in\N$ and $s_1,\dots,s_k\in\{1,\star\}$
and by using similar combinatorial arguments as in \cite{Thorbjornsen2000} one
can show that
\begin{equation}
\lim_{N\rightarrow\infty} \Var \tr (A_N^{s_1} \cdots
A_N^{s_k})=0.
\label{eq:zbieznoscb}
\end{equation}
Since $\limsup_{N\rightarrow\infty} \|T_N\|<\infty$ almost
surely \cite{Geman}, therefore there exists $R'>0$ such that for any integer
$k$ and $\epsilon>0$
\begin{equation}
\lim_{N\rightarrow\infty} P\big(\omega\in\Omega:
D_N+\sqrt{o} T_N(\omega)\in \Gamma_{R'}(x; k,N,\epsilon) \big) =1 .
\label{eq:prawieok}
\end{equation}
Furthermore since the convergence in (\ref{eq:zbieznosca}) and
(\ref{eq:zbieznoscb}) is uniform with respect to choice of the sequence
$(D_N)$, hence it is possible to find universal $R'$ for all choices of
$(D_N)$.

By comparing the densities of two measures on $\M_N^\sut$: the Lebesgue measure
$\vol^\sut$ and the distribution of the Gaussian random matrix $\sqrt{o} T_N$
we see that (\ref{eq:prawieok}) implies that for every $0<\delta<1$, every
integer $k$ and $\epsilon>0$ there exists $N_0$ such that for
$N>N_0$ we have that the volume of the set $\{m\in\M_N^\sut: D_N+m \in
\Gamma_{R'}(x;k,N,\epsilon) \}$ is bigger or equal to the volume of a
$N(N-1)$--dimensional ball with radius $\sqrt{(1-\delta)oN}$:
\begin{multline} \vol^\sut \{m\in\M_N^\sut: D_N+m \in
\Gamma_{R'}(x;k,N,\epsilon) \}\geq \\
\pi^{\frac{N(N-1)}{2}} \left[ \Gamma\left(
\frac{N (N-1)}{2}+1 \right)\right]^{-1} [(1-\delta)oN]^{\frac{N (N-1)}{2}}.
\label{eq:idenaobiad}
\end{multline}

Since the sequence $(D_N)$ was chosen arbitrarily, it follows that for every
$0<\delta<1$, every integer $k$ and $\epsilon>0$ there exists $N_0$ such that
for any $N>N_0$ and any $D_N\in\hat{\Gamma}^\dia_R(x; N, l_N,\theta_N)$
inequality (\ref{eq:idenaobiad}) holds.

Now it is enough to apply Proposition \ref{prop:pullback} to show that for
every $0<\delta<1$, every integer $k$ and $\epsilon>0$ we have
$$\chi_{R'}(x; k,\epsilon) \geq
\hat{\chi}^\dia_{R}(x)+\frac{1}{2}+\frac{\log 2\pi (1-\delta) o}{2}$$
what finishes the proof.
\end{proof}

\subsection{Comparison of inequalities for free entropy}
The inequality (\ref{eq:glowne}) contains a double integral, which is often
called logarithmic energy of a measure. A similar term appears in the formula
for the free entropy of a single selfadjoint operator which is due do
Voiculescu \cite{VoiculescuPart2}:
\begin{equation}
\chi^\sa(x)=\int_{\R} \int_{\R} \log |z_1-z_2| \dd\mu_x(z_1)
\dd\mu_x(z_2)+\frac{3}{4}+\frac{\log 2\pi}{2}.
\label{eq:entropiasamosprzezonego}
\end{equation}
It should be stressed that---despite the formal resemblance---the free
entropies in (\ref{eq:glowne}) and (\ref{eq:entropiasamosprzezonego}) are
different objects. Namely, in (\ref{eq:glowne}) we consider a free entropy
$\chi(x)$ of a non--selfadjoint random variable while in
(\ref{eq:entropiasamosprzezonego}) we consider a free entropy $\chi^\sa(x)$ of
a selfadjoint random variable which is defined by hermitian matrix
approximations.

On the other hand inequality (\ref{eq:glowne}) contains a term equal to the
logarithm of the offdiagonality, which can be regarded as a non--commutative
variance. A similar expression appears in the Voiculescu's inequality for a
free entropy of a non--selfadjoint variable \cite{VoiculescuPart2}:
\begin{equation}
\chi(x)\leq \log \left[ \pi e \big( \phi(|x|^2)-|\phi(x)|^2 \big)^2 \right].
\end{equation}

\section{Acknowledgements}
The research was conducted at Texas A\&{}M University on a
scholarship funded by Polish--US Fulbright Commission.
I acknowledge the support of Polish Research Committee
grant No.\ P03A05415.

\bibliography{biblio}

\end{document}